# Mathematical definition of general systems based on semiotics and set theory (part of statics)

Dandan Zou

Email: ddzou@mail.ustc.edu.cn

Abstract. According to the relation between objects and time, our category of general systems theory was divided into three parts: statics, kinematics, and dynamics. In this part, beginning with clarifying fundamental in epistemology and semiotics, we discussed the relationship among measurements, partitions and functions. Then we defined the terms of quantity and value by our understanding of functions. Farther more, the concept of relation quantity was coined to describe the relationship like distance and force. Finally, we present the structure which can afford the obligation of general system, and also the definitions of subsystem and isomorphism.

## 1. Introduction

Systems thinking can date back to antiquity. However, it is until 1920s that the concept of systems was first formally proposed by biologist Ludwig von Bertalanffy. Since then, the meaning of system has been getting more widely, and the application of the systems concept has become more widespread. In searching for providing a unified and formalized mathematical approach to all major systems concepts, Mesarovic and Takahara introduced a concept of general systems based on Cantor's set theory in the 1960s. Much effort for formal aspects of deterministic input-output systems, learning systems, decision-making systems, and goal-seeking systems has also been done recent years [1]. In this paper, we investigated only the static systems that our objects would not evolve with time. And in order to open the door towards kinematics and dynamics, the concept of quantity was proposed. To fully understand its meaning, rather than only restrict the definitions by words, we need to clarify how an individual perceives about the world.

## 2. Sensation and Object

All the cognitions of the external world derive from sensation, whereas other mental activities can only be the processes of the sensory material. Surprisingly, the sensation we perceived is not totally in disorder, but there exists something steady. This steady part transforms with time and measurement, usually exhibiting some kinds of stability, continuity, and the boundary with its environment sometimes may not be ideally closed and clear. Thereby from the sensation, we can extract various sensory objects, namely, entities. And the entities form into more complicated objects in our mind, after such process as abstraction, transplantation and so on.

## 3. Assignment and Sign

The human brain naturally has the ability of association which allows us assigning one object to another, as well as the assignment among many objects. Since such mechanism bases on the structure consists of neurons and their connections in our brain, assignment (namely correspondence) can not and need not be well defined. The operation of assignment which looks quite simple plays an important role during the rational thinking that follows the forming of objects. In a manner, rational knowledge can be considered as a kind of assignments.

The process of assigning something to other things is directed. We call the former as preimage (or inverse image), and latter as image. And undirected assignment can be considered as the combination of two inversed directed assignments. The assignment between preimage and image, which can be a representation of objective reality and also can have nothing to do with it, is arbitrary.

Sometimes it seems inconsistent that our cognitions of the world do not exactly come from our own sensation. For example, consider a thermometer in a room. We need not feeling to know the temperature in this room, but just use eyes to read its scale, or to be told by others. What we saw is not temperature, but rather a sign the thermometer conveyed to us. However, our own sensation indeed played an important role in such process, since without matched sense of touch the sign was meaningless for us.

From the discussion above, we can figure out two necessary conditions for any a sign. It must be a kind of assignment, and it must have been defined among particular people. In other words, a sign is an assignment defined by human being in a special group. Moreover, the image of the assignment (called signifier in semiotics) should have an advantage over preimage (signified) in the respect of conveying, although it is not necessary. In the following sections, we also call the image as a sign of preimage.

However there is still a question remained that how we can define a sign or symbol. Usually, the preimage and image of a sign is assigned by using a sequence of different signs. It can be proved that the above method which is called symbolic definition could not be the only method to define signs or symbols, since there must be circular definitions if all signs or symbols can be defined in this way. And nonsymbolic definition (the definition without using other symbols) can emerge between two individuals with shared experience.

## 4. Measurement and Function

Measurement is essentially a kind of set partition [2, 3]. By following the same operation, we can compare and distinguish two different objects in a set. For example, we can distinguish between different colors of objects, and assign them to corresponding signs of red, blue, orange and so on. Similarly, we can distinguish the distances of two geometrical objects, and assign them to corresponding real numbers which can be an image of a sign. Comparing one object with another can be considered as a process of measurement, no matter which sense we used or which method we operated. More generally, even the algorithm of calculating parity, which divides the set of all positive integers into two subsets, is some sort of measurement without any sensation or behavior.

For specifying our realm in this paper, we are talking about only static measurement that objects will not change and measurement can not affect objects in a classical situation.

Definition  [Function and Operation] For A and R two nonempty sets and k a nonnegative integer we define $A^0 = \{\emptyset\}$, and, for k> 0, $A^K$ is the set of k-tuples of elements from A. According to a defined rule F of correspondence, if every element of $A^K$ is associated with exactly one element in R, then F is said to be a k-ary function from A to R, written F: $A^K \to R$. If R=A, then F is said to be an n-ary operation on A [4].

Theorem  For any a function F: $A^K \to R$, F can induce a partition of set $A^K$. [5]

This theorem tells us that a function can give expression to a kind of partition or measurement. And if only the element in the range is different from domain (in fact, even this is not a rigorous condition), it can be chose arbitrary. It is easy to find examples that the values of function only play a role of signs in both mathematics and everyday life.

However, we can perceive colors without knowing any language. Partition and assigning the subsets to signs can be two separate processes. And the latter may not exist in some measurement phenomena.

## 5. Quantity and Value [6]

A quantity is a rule of correspondence established between the set of objects and the set of its signs. For a function F: $A^K \to R$, which is specially used for expressing only one kind of partition, R is essentially a set of signs and can denote this particular function.

Definition  Suppose A is a nonempty set of (distinct) objects, and R is a nonempty set. If $A \cap R = \emptyset$ and there is a function F: $A^K \to R$, then R is said to be a quantity of A, while the elements in R are called as values of A in R, written $A^K \to R$.

For k=1, R is said to be an attribute of A, and elements in R are called as characters of A in R;

For k>1, R is said to be a relation quantity, and elements in R are called as relation values (or relations) of A in R.

If F is a bijection, then R is said to be a tag of A, written $A^K \leftrightarrow R$.

Suppose $A^K \to R, x_1, x_2, \cdots, x_k \in A, R_i \in R$. if a tuple $<x_1, x_2, \cdots, x_k>$ is assigned to $R_i$, then it is denoted by $R(x_1, x_2, \cdots, x_f) = R_i$, or $R_i(x_1, x_2, \cdots, x_f)$;

Suppose $A^K \to L, L_i \in L$, then a subset of $A^K$ is assigned to $L_i$. The subset can be denoted by $\{x \in A^K | L(x) = L_i\}$, or $\{A^K | L_i\}$.

Some relationship in our everyday life like marriage (between men and women) or belonging (between elements and sets), can be described by the concept relation in traditional set theory. However, it is our intention to coin the concept of relation quantity that there is also some relationship like distance in geometry or force in mechanics that can not be described by such customary way. And as will be discussed, the relation in set theory can also be contained in the relation quantity by considering it as a two-element quantity.

Theorem    If $B \subset A$, and $A^n \to L$, then $B^n \to L$.

Theorem    Suppose n, m are positive integers and n<m. If $A^n \to L$, then $A^m \to L$.

## 6. Quantities

Correspondences that established between two quantities on the mutual set of objects can be divided into three types.

Definition    A is a nonempty set of objects. For L and S two n-ary quantities of A, $L_i$ and $S_i$ are elements of L and S, respectively. If $\{A^n | L_i\} \cap \{A^n | S_i\} \neq \emptyset$, then there is a Z-Z type correspondence between $L_i$ and $S_i$, written $L_i |\xleftrightarrow{A}| S_i$, otherwise written $L_i * \xleftrightarrow{A} * S_i$.

If $\{A^n | L_i\} \neq \emptyset$, and $\{A^n | L_i\} \subset \{A^n | S_i\}$, then there is a Z->Z type correspondence between $L_i$ and $S_i$, written $L_i \xrightarrow{A} S_i$.

Moreover if $L_i \xrightarrow{A} S_i$, and $S_i \xrightarrow{A} L_i$, then there is a Z<->Z type correspondence between $L_i$ and $S_i$, written $L_i \xleftrightarrow{A} S_i$.

Theorem    If $L_i \xrightarrow{A} S_i$, then $L_i |\xleftrightarrow{A}| S_i$.

**Definition** If $A^n \rightarrow L$, then L is said to be a constant quantity on A while $(\exists L_i \in L)(\forall x \in A^n)(L(x) = L_i)$.

**Definition** If $(\forall L_i \in L, \forall S_i \in S)(L_i |\xleftrightarrow{A}| S_i)$, then L is independent of S on A, written $L * \xleftrightarrow{A} * S$.

**Example** In analytic geometry, all points can be the set of objects and any a coordinate axis can be a quantity. It is obvious that any two axes are independent of each other.

S is a function (or dependent) of L on A, written $L \xrightarrow{A} S$ : $A^n \rightarrow L$, $A^n \rightarrow S$, $(\forall L_i \in L)(\exists S_i \in S)(\text{if}\{A^n | L_i\} \neq \varnothing$, then $L_i \xrightarrow{A} S_i)$.

If $L \xrightarrow{A} S$, and $S \xrightarrow{A} L$, then L and S are equivalent on A, written $L \xleftrightarrow{A} S$.

It is an equivalence relation on the set of all quantities on A, which means that L and S represent the same partition on A.

The set of all partitions of A is denoted by $\Pi(A)$. The set of all quantities on A is denoted by Qu(A).

A partial order $\leq$ can be defined on Qu(A) by L$\leq$S iff $L \xrightarrow{A} S$, as well as $\leq$ defined on $\Pi(A)$ by $\pi_1 \leq \pi_2$ iff each block of $\pi_1$ is contained in some block of $\pi_2$ [4].

**Definition** Let P be a subset of Qu(A). A element U in Qu(A) is an upper bound of P if p$\leq$U for every p in P. A element U in Qu(A) is the least upper bound of P or supremum of P if U is an upper bound of P, and p$\leq$m for every p in P implies U$\leq$m. Similarly we can define lower bound and the greatest lower bound.

**Definition** $A^n \rightarrow L$, $A^n \rightarrow S$, $L \wedge S$ is the greatest lower bound of L and S iff it such that:

1. $(L \wedge S) \xrightarrow{A} L$; 2. $(L \wedge S) \xrightarrow{A} S$;
3. $(\forall k)(A \text{-}> k, \text{if } k \xrightarrow{A} L, k \xrightarrow{A} S, \text{then } k \xrightarrow{A} (L \wedge S))$.

Similarly $L \vee S$ can be defined.

Definition    If $A^n \leftrightarrow L \wedge S \wedge J$, and L, S, J are not constant quantities, and L, S, J are independent of each other, then {L, S, J} is said to be a n-ary complete set of A.

Theorem    If {L} is an n-ary complete set of A, then there is not an n-ary quantity S such that {L, S} is an n-ary complete set of A.

Theorem    If {L, S, J} is an n-ary complete set of A, then $L \wedge S \wedge J \xrightarrow{A} U$ for any n-ary quantity U on A.

The general situation of a complete set $\{J_1, J_2, \cdots J_m\}$ can be discussed similarly.

The dividing line of objects and signs is not absolute. As a set of signs, a quantity can also be consider as a set of objects and have its own quantities. For example, there are many quantities on the set of natural numbers, which can be used as a quantity of a set of objects, such as total order and arithmetic operations.

Theorem    If $A^n \to L$, and $L^m \to R$, then R can be an n×m-ary quantity on A.

## 7. System

Definition    Suppose V is a nonempty set and $\vec{R}$ is a tuple of which elements are consist of quantities defined on V. Then the ordered pair $<V, \vec{R}>$ is said to be a system, written SYS$<V, \vec{R}>$. V is said to be a vertex set of the system, written V_SYS. $\vec{R}$ is said to be a measure set of the system, written M_SYS. V and $\vec{R}$ are called as structure of the system. More strictly, $\vec{R}$ should contain at least one relation quantity. And sometimes, partition can also be considered into the measure set.

Axiom    Any a system is an object that can be an element of another set or system.

Suppose SA and SB is two different systems. If SA is an element of V_SB, then SA is a subsystem of SB.

If $V\_SA \subset V\_SB$, and $M\_SA = M\_SB$, then SA is part of SB.

If $V\_SA = V\_SB$, and $M\_SA \subset M\_SB$, then SB is concretion of SA, or SA is abstraction of SB.

Let R be a quantity on SB, and SA be part of SB. If R is a constant quantity on SA, then SA is uniform about R.

Definition  Two systems SA and SB are isomorphic if there is a bijection F from SA to SB such that for every x in $V\_SA$ and every quantity R in $M\_SA$ the following equation holds: $R(x)=F[F[R](F[x])]$. Such F is called isomorphism. If only F: $V\_SA \rightarrow V\_SB$ is not one-to-one, then F: SA$\rightarrow$SB is a homomorphism from SA to SB.

Definition  A is a set of systems, and $A \rightarrow R$. $\forall S_i \in A$, if R is a function related to structure of $S_i$, then R is said to be a structure attribute of A.

Any an attribute of physical objects is always determined by its structure, since under the ideal measurement situation there is no reason of measuring results to be different for the identical inside of structures and outside of measuring conditions or methods. However, it is only a conjecture without more strict proof presented here. And for objects in our mind, sometimes they really can be assigned to any value freely without considering their structures inside.

## 8. Future work: kinetic and dynamical systems

The dynamical systems, have been  riginally by Poincare. The property of input and output for systems will be based on it.

**Acknowledgment**

I thank K.L.Wang for thoughtful comments, and Q.Y.Mo, L.Sun for reading and rectifying grammatical errors in an early draft.